\documentclass{amsart}

\usepackage{amssymb}
\usepackage{latexsym}
\usepackage{amsmath}
\usepackage{euscript}
\usepackage{hyperref}

 \newtheorem{thm}{Theorem}[section]
 \newtheorem{cor}[thm]{Corollary}
 \newtheorem{lem}[thm]{Lemma}
 \newtheorem{prop}[thm]{Proposition}
 \theoremstyle{definition}
 \newtheorem{defn}[thm]{Definition}
 \theoremstyle{remark}
 \newtheorem{rem}[thm]{Remark}
 
 \numberwithin{equation}{section}

\def\sS{{\mathfrak S}}

      \def\dC{{\mathbb C}}

   \def\dN{{\mathbb N}}   
      \def\dR{{\mathbb R}}

   \def\cB{{\mathcal B}}   
      
\def\cG{{\mathcal G}}   \def\cH{{\mathcal H}}   
   \def\cK{{\mathcal K}}

\newcommand{\dom}{{\mathrm{dom\,}}}
\newcommand{\ran}{{\mathrm{ran\,}}}

\DeclareMathOperator{\supp}{supp}

\newcommand{\Lar}{A_{\alpha}}
\newcommand{\AN}{A_{\rm N}}

\DeclareMathOperator\dist{dist}

\newcommand\eps{\varepsilon}

\newcommand\ov{\overline}
\newcommand\wt{\widetilde}

\newcommand{\defeq}{\mathrel{\mathop:}=}

\begin{document}
%
%
%
%
%
%
%
%
%
\title[Schatten-von Neumann estimates for resolvent differences]{Schatten-von Neumann estimates for resolvent differences of Robin Laplacians on a half-space}
\author[V.~Lotoreichik \and J.~Rohleder]{Vladimir Lotoreichik \and Jonathan Rohleder}

\address{Technische Universit\"at Graz\\
Institut f\"ur Numerische Mathematik\\
Steyrergasse 30\\
8010 Graz, Austria}

\email{rohleder@math.tugraz.at \and lotoreichik@math.tugraz.at}




\begin{abstract}
The difference of the resolvents of two Laplacians on a half-space subject to Robin boundary conditions is studied. In general this difference is not compact, but it will be shown that it is compact and even belongs to some Schatten-von-Neumann class, if the coefficients in the boundary condition are sufficiently close to each other in a proper sense. In certain cases the resolvent difference is shown to belong even to the same Schatten-von Neumann class as it is known for the resolvent difference of two Robin Laplacians on a domain with a compact boundary. 
\end{abstract}

\maketitle

\section{Introduction}

Schatten-von Neumann properties for resolvent differences of elliptic operators on domains have been studied basically since M.\,Sh.~Birman's famous paper~\cite{B62}, which appeared fifty years ago and was followed by important contributions of M.\,Sh.~Birman and M.\,Z.~Solomjak as well as of G.~Grubb, see~\cite{BS80, G84a, G84}; moreover, the topic has attracted new interest very recently, see~\cite{BLLLP10,BLL10,G11,G11-2, M10}. Recall that a compact operator belongs to the Schatten-von Neumann class $\sS_p$ (weak Schatten-von Neumann class $\sS_{p, \infty}$) of order $p > 0$ if the sequence of its singular values is $p$-summable (is $O (k^{- 1/p})$ as $k \to \infty$); see Section~\ref{sec:main} for more details. The objective of the present paper is to study the resolvent difference of two (in general non-selfadjoint) Robin Laplacians on the half-space $\dR^{n + 1}_+ = \{ (x',x_{n+1})^{\rm T}\colon x'\in\dR^n, x_{n+1} > 0\}$, $n \geq 1$, of the form
\begin{equation}\label{eq:dom}
	 \Lar f = - \Delta f, \quad \dom(\Lar) = 
	\left\{ f\in H^2 (\dR^{n + 1}_+) \colon \partial_\nu f|_{\dR^n}
	= \alpha f|_{\dR^n} \right\},
\end{equation}
in $L^2 (\dR^{n + 1}_+)$ with a function $\alpha : \dR^n \to \dC$ belonging to the Sobolev space $W^{1, \infty} (\dR^n)$, i.e., $\alpha$ is bounded and has bounded partial derivatives of first order; here $f |_{\dR^n}$ is the trace of a function $f$ at the boundary $\dR^n$ of $\dR^{n + 1}_+$ and $\partial_\nu f |_{\dR^n}$ is the trace of the normal derivative of $f$ with the normal pointing outwards of $\dR^{n+1}_+$. We emphasize that, as a special case, our discussion contains the resolvent difference of the selfadjoint operator with a Neumann boundary condition and a Robin Laplacian. If the half-space $\dR^{n + 1}_+$ is replaced by a domain with a compact, smooth boundary, it is known that for real-valued $\alpha_1$ and $\alpha_2$ the operators $A_{\alpha_1}$ and $A_{\alpha_2}$ are selfadjoint and the difference of their resolvents
\begin{align}\label{resdiff}
 (A_{\alpha_2} - \lambda)^{-1} - (A_{\alpha_1} - \lambda)^{-1}
\end{align}
belongs to the class $\sS_{\frac{n}{3}, \infty}$; see~\cite{BLL10} and~\cite{BLLLP10,G11-2}, where also more general elliptic differential expressions and certain non-selfadjoint cases are discussed. 

On the half-space $\dR^{n + 1}_+$ the situation is fundamentally different. Here, in general, the resolvent difference~\eqref{resdiff} is not even compact. For instance, if $\alpha_1 \neq  \alpha_2$ are real, positive constants, the essential spectra of $A_{\alpha_1}$ and $A_{\alpha_2}$ are given by $[- \alpha_1^2, \infty)$ and $[- \alpha_2^2, \infty)$, respectively. Consequently, in this case the difference~\eqref{resdiff} cannot be compact. Nevertheless, the main results of the present paper show that under the assumption of a certain decay of the difference $\alpha_2 (x) - \alpha_1 (x)$ for $|x| \to \infty$, compactness of the resolvent difference in~\eqref{resdiff} can be guaranteed, and that this difference belongs to $\sS_p$ or $\sS_{p, \infty}$ for certain $p$, if $\alpha_2 - \alpha_1$ has a compact support or belongs to $L^q (\dR^n)$ for some $q$. It is a question of special interest under which assumptions on $\alpha_2-\alpha_1$  the difference~\eqref{resdiff} belongs to $\sS_{\frac{n}{3},\infty}$, i.e., to the same class as in the case of a domain with a compact boundary. Our results show that if $\alpha_2 - \alpha_1$ has a compact support this is always true, and that in dimensions $n > 3$ a sufficient condition is
\begin{equation*}
 \alpha_2 -\alpha_1\in L^{n/3}(\dR^n).
\end{equation*}
If $\alpha_2 - \alpha_1 \in L^{p}(\dR^n)$ with $p \ge 1$ and $ p > n/3$, we show that the resolvent difference in \eqref{resdiff} belongs to the larger class $\sS_p \supsetneq \sS_{\frac{n}{3},\infty}$. In dimensions $n = 1,2$ for $\alpha_2 - \alpha_1 \in L^1 (\dR^n)$ the difference~\eqref{resdiff} belongs to the trace class $\sS_1$. As an immediate consequence, for $n=1,2$ and real-valued $\alpha_1$, $\alpha_2$ with $\alpha_2 - \alpha_1 \in L^1 (\dR^n)$ wave operators for the pair $\{A_{\alpha_1},A_{\alpha_2}\}$ exist and are complete, which is of importance in scattering theory. 
Two further corollaries of our results concern the case that $A_{\alpha_1}$ is the Neumann operator, i.e., $\alpha_1 = 0$: on the one hand, if $\alpha_2$ is real-valued, under our assumptions the Neumann operator and the absolutely continuous part of $A_{\alpha_2}$ are unitarily equivalent, cf.~\cite{MN11}; on the other hand, with the help of recent results from perturbation theory for non-selfadjoint operators, see~\cite{DHK09,H11}, we conclude some statements on the accumulation of the (in general non-real) eigenvalues in the discrete spectrum of $A_{\alpha_2}$. 

Our results complement and extend the result by M.\,Sh. Birman in~\cite{B62}. He considers a realization of a symmetric second-order elliptic differential expression on an unbounded domain with combined boundary conditions, a Robin boundary condition on a compact part and a Dirichlet boundary condition on the remaining non-compact part of the boundary, and showed that the resolvent difference of  the described realization and  the realization with a Dirichlet boundary condition on the whole boundary belongs to the class $\sS_{\frac{n}{2},\infty}$. It is remarkable that in our situation in some cases the singular values converge faster than in the situation Birman considers. This phenomenon is already known for domains with compact boundaries, when a Neumann boundary condition instead of a Dirichlet boundary condition is considered; see~\cite{BLLLP10}.

It is worth mentioning that all results in this paper on compactness and Schatten-von Neumann estimates remain valid for $- \Delta$ replaced by a Schr\"odinger differential expression $- \Delta + V$ with a real-valued, bounded potential $V$ and the proofs are completely analogous.

Our considerations are based on an abstract concept from the extension theory of symmetric operators, namely, the notion of quasi-boundary triples, which was introduced by J.~Behrndt and M.~Langer in~\cite{BL07} and has been developed further by them together with the first author of the present paper in~\cite{BLL10}. The key tool provided by the theory of quasi-boundary triples is a convenient factorization of the resolvent difference in \eqref{resdiff}. For the proof of our main theorem we combine this factorization with results on the compactness of the embedding of $H^1 (\Omega)$ into $L^2(\Omega)$ for $\Omega$ being a (possibly unbounded) domain of finite measure and with $\sS_p$- and $\sS_{p,\infty}$-properties of the operator $\sqrt{|\alpha_2-\alpha_1|}(I - \Delta_{\dR^n})^{-3/4}$; the proof of the most optimal $\sS_{\frac{n}{3},\infty}$-estimate is based on an asymptotic result proved by M.\ Cwikel in~\cite{C77}, conjectured earlier by B.\ Simon in \cite{S76}.

A short outline of this paper looks as follows. In Section~\ref{sec:prelim} we give an overview of some known results on quasi-boundary triples which are used in the further analysis and provide a quasi-boundary triple for the Laplacian on the half-space; furthermore we prove that for each two coefficients $\alpha_1, \alpha_2$ the operators $A_{\alpha_1}$ and $A_{\alpha_2}$ have joint points in their resolvent sets and that $A_\alpha$ is selfadjoint if and only if $\alpha$ is real-valued. In Section~\ref{sec:main} we establish sufficient conditions for the resolvent difference~\eqref{resdiff} to be compact or even to belong to certain Schatten-von Neumann classes. The paper concludes with some corollaries of the main results.

Let us fix some notation.  If $T$ is a linear operator  from a Hilbert space $\cH$ into a Hilbert space $\cG$ we denote by $\dom T$, $\ran T$, and $\ker T$ its domain, range, and kernel, respectively. If $T$ is densely defined, we write $T^*$ for the adjoint operator of $T$. If $\Theta$ and  $\Lambda$ are linear relations from $\cH$ to $\cG$, i.e., linear subspaces of $\cH \times \cG$, we define their sum to be
\begin{align*}
 \Theta + \Lambda = \big\{ \{f, g_\Theta + g_\Lambda\} : \{f, g_\Theta\} \in \Theta, \{f, g_\Lambda\} \in \Lambda \big\}.
\end{align*}
We write $T\in\cB(\cH,\cG)$, if $T$ is a bounded, everywhere defined operator from $\cH$ into $\cG$; if $\cG = \cH$ we simply write $T \in \cB (\cG)$. For a closed operator $T$ in $\cH$ we denote by $\rho(T)$ and $\sigma(T)$ its resolvent set and spectrum, respectively. Moreover, $\sigma_{\rm d} (T)$ denotes the discrete spectrum of $T$, i.e., the set of all eigenvalues of $T$ which are isolated in $\sigma (T)$ and have finite algebraic multiplicity, and $\sigma_{\rm ess} (T)$ is the essential spectrum of $T$, which consists of all points $\lambda \in \dC$ such that $T - \lambda$ is not a semi-Fredholm operator. Finally, for a bounded, measurable function $\alpha : \dR^n \to \dC$ we denote its norm by $\| \alpha \|_\infty = \sup_{x \in \dR^n} | \alpha (x)|$. Furthermore, for simplicity we identify $\alpha$ with the corresponding multiplication operator in $L^2 (\dR^n)$.

\section{Quasi-boundary triples and Robin Laplacians on a half-space}
\label{sec:prelim}

In this section we provide some general facts on quasi-boundary triples as introduced in~\cite{BL07}. Afterwards we apply the theory to the Robin Laplacian in~\eqref{eq:dom}. Let us start with the basic definition. 

\begin{defn}
\label{def:qbt}
Let $A$ be a closed, densely defined, symmetric operator in a Hilbert space $(\cH,(\cdot,\cdot)_\cH)$.
We say that $\{\cG,\Gamma_0,\Gamma_1\}$ is a \emph{quasi-boundary triple}
for $A^*$, if $T\subset A^*$ is an operator satisfying $\ov T  = A^*$ and  $\Gamma_0$ and $\Gamma_1$ are linear mappings defined on
$\dom T$ with values in the Hilbert space $(\cG,(\cdot,\cdot)_\cG)$ such that the following conditions are satisfied.
\begin{enumerate}
 \item $\Gamma \defeq \binom{\Gamma_0}{\Gamma_1}\colon
\dom T\rightarrow\cG\times\cG$ has a dense range.
 \item The {\em abstract Green identity}
\begin{equation}\label{green1}
  (Tf,g)_{\cH}-(f,Tg)_{\cH}
  =(\Gamma_1 f,\Gamma_0 g)_{\cG}-(\Gamma_0 f,\Gamma_1 g)_{\cG}
\end{equation}
holds for all $f,g\in \dom T$.
 \item $A_0 := T \upharpoonright \ker\Gamma_0 = A^*\upharpoonright\ker\Gamma_0$ is selfadjoint.
\end{enumerate}
\end{defn}

We set $\cG_i = \ran \Gamma_i$, $i = 0, 1$. Note that the definition of a quasi-boundary triple as given above is only a special case of the original one given in~\cite{BL07} for the adjoint of a closed, symmetric linear relation $A$. We remark that if $\{\cG,\Gamma_0,\Gamma_1\}$ is a quasi-boundary triple with the additional property $\cG_0 =\cG$, then $\{\cG,\Gamma_0,\Gamma_1\}$ is a generalized boundary triple in the sense of~\cite{DM95}. Let us also mention that a quasi-boundary triple for $A^*$ exists if and only if
the deficiency indices $\dim\ker(A^*\mp i)$ of $A$ coincide. If $\{\cG, \Gamma_0, \Gamma_1 \}$ is a quasi-boundary triple for $A^*$ with $T$ as in Definition~\ref{def:qbt}, then $A$ coincides with $T \upharpoonright {\ker\Gamma}$.

The next proposition contains a sufficient
condition for a
triple $\{\cG,\Gamma_0,\Gamma_1\}$ to be a quasi-boundary triple. For a proof see~\cite[Theorem~2.3]{BL07}; cf. also \cite[Theorem 2.3]{BL11}.
 
\begin{prop}\label{suff_cond_qbt}
Let $\cH$ and $\cG$ be Hilbert spaces and let $T$ be a linear operator in $\cH$.
Assume that
$\Gamma_0,\Gamma_1\colon \dom T\rightarrow\cG$ are linear mappings such that the
following
conditions are satisfied:
\begin{itemize}\setlength{\itemsep}{1.2ex}
\item[(a)]
$\Gamma\defeq \binom{\Gamma_0}{\Gamma_1}\colon \dom T\rightarrow\cG\times\cG$
has a dense range.
\item[(b)]
The identity \eqref{green1} holds for all $ f,g\in \dom T$.
\item[(c)] $T\upharpoonright {\ker\Gamma_0}$ contains a selfadjoint operator.
\end{itemize}
Then 
$A \defeq T \upharpoonright {\ker\Gamma}$ is a closed, densely defined, symmetric operator in $\cH$ and
$\{\cG,\Gamma_0,\Gamma_1\}$
is a quasi-boundary triple for $A^*$.
\end{prop}

Let us recall next the definition of two related analytic objects, the $\gamma$-field and the Weyl function associated with a quasi-boundary triple.

\begin{defn}
Let $A$ be a closed, densely defined, symmetric operator in $\cH$ and let
$\{\cG,\Gamma_0,\Gamma_1\}$ be a quasi-boundary triple for $A^*$ with $T$ as in Definition~\ref{def:qbt} and $A_0 = T \upharpoonright {\ker\Gamma_0}$.
Then the operator-valued functions $\gamma$ and $M$ defined by
\begin{equation*}
  \gamma(\lambda) \defeq \bigl(\Gamma_0\upharpoonright\ker(T-\lambda)\bigr)^{-1}
  \,\,\,\text{and}\,\,\,
  M(\lambda) \defeq \Gamma_1\gamma(\lambda),\quad
  \lambda\in\rho(A_0),
\end{equation*}
are called the $\gamma$\emph{-field} and the \emph{Weyl function}, respectively, corresponding to
the triple $\{\cG,\Gamma_0,\Gamma_1\}$.
\end{defn}

These definitions coincide with the definitions of the $\gamma$-field and
the Weyl function in the case that $\{\cG,\Gamma_0,\Gamma_1\}$ is an ordinary 
boundary triple, see~\cite{DM91}. It is an immediate consequence of the decomposition
\begin{align}\label{decompo}
 \dom T = \dom A_0 \dotplus \ker (T - \lambda) = \ker \Gamma_0 \dotplus \ker (T - \lambda), \quad \lambda \in \rho(A_0),
\end{align}
that the mappings $\gamma (\lambda)$ and $M (\lambda)$ are well-defined. Note that for each $\lambda \in \rho(A_0)$, $\gamma(\lambda)$ maps $\cG_0$ onto $\ker (T - \lambda) \subset \cH$ and $M(\lambda)$ maps $\cG_0$ into $\cG_1$. Furthermore, it follows immediately from the definitions of $\gamma (\lambda)$ and $M(\lambda)$ that
\begin{align}\label{NDpropertyAbstr}
 \gamma (\lambda) \Gamma_0 f_\lambda = f_\lambda \quad \text{and} \quad M(\lambda) \Gamma_0 f_\lambda = \Gamma_1 f_\lambda, \quad f_\lambda \in \ker (T - \lambda),
\end{align}
holds for all $\lambda \in \rho(A_0)$.

In the next proposition we collect some properties of the $\gamma$-field
and the Weyl function; all statements can be found in \cite[Proposition~2.6]{BL07}.

\begin{prop} \label{gammaprop}
Let $A$ be a closed, densely defined, symmetric operator in a Hilbert space $\cH$ and let
$\{\cG,\Gamma_0,\Gamma_1\}$ be a quasi-boundary triple for $A^*$ with
$\gamma$-field $\gamma$ and Weyl function $M$. Denote by $A_0$ the restriction of $A^*$ to $\ker \Gamma_0$. Then for $\lambda \in \rho (A_0)$ the following assertions hold.
\begin{itemize}\setlength{\itemsep}{1.2ex}
\item[(i)] 
$\gamma(\lambda)$ is a bounded, densely defined operator from $\cG$ into $\cH$.
\item[(ii)] 
The adjoint of $\gamma(\overline\lambda)$ can be expressed as 
\begin{equation*}
\gamma(\ov\lambda)^* = \Gamma_1(A_0-\lambda)^{-1} \in\cB(\cH,\cG).
\end{equation*}
\item[(iii)] 
$M(\lambda)$ is a densely defined, in general unbounded operator in $\cG$, whose range is contained in $\cG_1$ and which satisfies $M(\ov \lambda) \subset M(\lambda)^*$.
\end{itemize}
\end{prop}

A quasi-boundary triple provides a parametrization for a class of extensions of a closed, densely defined, symmetric operator $A$. If $\{\cG, \Gamma_0, \Gamma_1\}$ is a quasi-boundary triple for $A^*$ with $T$ as in Definition~\ref{def:qbt} and $\Theta$ is a linear relation in $\cG$, we denote by $A_\Theta$ the restriction of $T$ given by
\begin{equation}\label{athetaop}
  A_\Theta f = T f, \quad \dom A_\Theta = \left\{ f\in \dom T \colon \bigl(\begin{smallmatrix}\Gamma_0 f\\ \Gamma_1 f\end{smallmatrix}\bigr)\in\Theta \right\}.
\end{equation}
In contrast to the case of an ordinary boundary triple, this parametrization does not cover all extensions of $A$ which are contained in $A^*$, and selfadjointness of  $\Theta$
does not imply selfadjointness or essential selfadjointness of $A_\Theta$;
cf.\ \cite[Proposition~4.11]{BL07} for a counterexample and
\cite[Proposition~2.4]{BL07}. The following proposition shows that under certain conditions $\rho (A_\Theta)$ is non-empty, and it implies a sufficient condition for selfadjointness of $A_\Theta$. It also provides a formula of Krein type for the resolvent difference of $A_\Theta$ and $A_0$. In the present form the proposition is a special case of~\cite[Theorem~2.8]{BL07}.

\begin{prop} \label{prop.krein}
Let $A$ be a closed, densely defined, symmetric operator in $\cH$ and let $\{\cG,\Gamma_0,\Gamma_1\}$
be a quasi-boundary triple for $A^*$ with $A_0 = A^* \upharpoonright\ker\Gamma_0$. Let $\gamma$ be the corresponding $\gamma$-field and $M$ the corresponding Weyl function.
Furthermore, let $\Theta$ be a linear relation in $\cG$ and assume that $(\Theta - M (\lambda))^{-1} \in \cB (\cG)$ is satisfied for some $\lambda \in \rho(A_0)$. Then $\lambda \in \rho (A_\Theta)$ and 
\begin{equation*} 
  (A_\Theta-\lambda)^{-1} - (A_0-\lambda)^{-1}
  = \gamma(\lambda)\bigl(\Theta-M(\lambda)\bigr)^{-1}\gamma(\ov\lambda)^*
\end{equation*}
holds.
\end{prop}

In order to construct a specific quasi-boundary triple for $- \Delta$ on the half-space $\dR^{n + 1}_+$, $n \geq 1$, let us recall some basic facts on traces of functions from Sobolev spaces. For proofs and further details see, e.g.,~\cite{AF03,G09,LM72}. We denote by $H^s (\dR^{n + 1}_+)$ and $H^s (\dR^n)$ the $L^2$-based Sobolev spaces of order $s \geq 0$ on $\dR^{n + 1}_+$ and its boundary $\dR^n$, respectively. The closure in $H^s (\dR^{n + 1}_+)$ of the space of infinitely-differentiable functions with a compact support is denoted by $H_0^s (\dR^{n + 1}_+)$. The trace map \linebreak $C^\infty (\overline{\dR^{n+1}_+}) \ni f \mapsto f|_{\dR^n} \in C^\infty (\dR^n)$ and the trace of the derivative $C^\infty (\overline{\dR^{n+1}_+}) \ni f \mapsto \partial_\nu f|_{\dR^n} = - \frac{\partial f}{\partial x_{n + 1}}\big |_{\dR^n} \in C^\infty (\dR^n)$ in the direction of the normal vector field pointing outwards of $\dR^{n+1}_+$ extend by continuity to $H^s(\dR^{n+1}_+), s > 3/2$, such that the mapping
\begin{equation*}
H^s(\dR_+^{n+1}) \ni f \mapsto \binom{f|_{\dR^n}}{\partial_\nu f|_{\dR^n}} \in H^{s-1/2}(\dR^n)\times H^{s-3/2}(\dR^n)
\end{equation*}
is well-defined and surjective onto $H^{s-1/2}(\dR^n)\times H^{s-3/2}(\dR^n)$. Moreover, this mapping can be extended to the spaces
\[
H^s_\Delta(\dR^{n+1}_+) := \bigl\{ f\in H^s(\dR^{n+1}_+) \colon \Delta f\in L^2(\dR^{n+1}_+)\bigr\}, \quad s \ge 0;
\]
cf.~\cite{Fr67, G08, G09, LM72}. We remark that for $s \ge 2$ the latter space coincides with the usual Sobolev space $H^s (\dR^n)$. In contrast to the case $s \geq 2$, the mapping
\begin{equation*}
H^s_\Delta(\dR_+^{n+1}) \ni f \mapsto \binom{f|_{\dR^n}}{\partial_\nu f|_{\dR^n}} \in H^{s-1/2}(\dR^n)\times H^{s-3/2}(\dR^n), \quad s\in [0,2),
\end{equation*}
is not surjective onto the product $H^{s-1/2}(\dR^n)\times H^{s-3/2}(\dR^n)$, but the separate mappings
\begin{align*}
&H^s_\Delta(\dR_+^{n+1}) \ni f \mapsto f|_{\dR^n}\in H^{s-1/2}(\dR^n), \quad s\in [0,2),
\end{align*}
and
\begin{align}\label{NTrace}
H^s_\Delta(\dR_+^{n+1}) \ni f\mapsto  \partial_\nu f|_{\dR^n}\in H^{s-3/2}(\dR^n),\quad s\in [0,2),
\end{align}
are surjective onto $H^{s-1/2}(\dR^n)$ and $H^{s-3/2}(\dR^n)$, respectively.

Let us introduce the operator realizations of $- \Delta$ in $L^2(\dR^{n+1}_+)$ given by
\begin{equation}
\label{eq:A}
	 A f = -\Delta f,\quad 
	\dom A = H^2_0(\dR^{n+1}_+),
\end{equation}
and
\begin{equation*}
	Tf = -\Delta f,\quad 	\dom T = H^{3/2}_{\Delta}(\dR^{n+1}_+), 
\end{equation*}
and the boundary mappings $\Gamma_0$ and $\Gamma_1$ defined by
\begin{equation}
\label{eq:bmap}
\Gamma_0f = \partial_\nu f|_{\dR^n}, \quad 
\Gamma_1f = f|_{\dR^n},\qquad f\in\dom T.
\end{equation}
Furthermore, let us mention that the {\em Neumann operator}
\begin{align}\label{eq:DirNeu}
	 A_{\rm N} f = -\Delta f,\quad \dom A_{\rm N}  = \big\{f\in H^2(\dR^{n+1}_+)\colon 
	 \partial_\nu f \big|_{\dR^n} = 0 \big\},
\end{align}
is selfadjoint in $L^2(\dR^{n + 1}_+)$ and its spectrum is given by $\sigma (A_{\rm N}) = [0, \infty)$; see, e.g.,~\cite[Chapter~9]{G09}. We prove now that the mappings $\Gamma_0$ and $\Gamma_1$ in~\eqref{eq:bmap} provide a quasi-boundary triple for the operator $A^*$ with $A_0 := A^* \upharpoonright \ker \Gamma_0 = A_{\rm N}$.

\begin{prop}\label{qbt}
The operator $A$ in \eqref{eq:A} is closed, densely defined, and symmetric, and the triple $\{\cG,\Gamma_0,\Gamma_1\}$ with  $\cG = L^2(\dR^n)$ and $\Gamma_0,\Gamma_1$ defined in~\eqref{eq:bmap} is a quasi-boundary triple for $A^*$.  Moreover, $A^* \upharpoonright \ker \Gamma_0 = A_{\rm N}$ holds. For $\lambda \in \rho(A_{\rm N})$ the associated $\gamma$-field is given by the {\em Poisson operator}
\begin{align}\label{Poisson}
	 \gamma (\lambda) \partial_\nu f_\lambda|_{\dR^n} = 
	 f_\lambda, \quad f_\lambda \in \ker (T - \lambda), 
\end{align}
and the associated Weyl function is given by the {\em Neumann-to-Dirichlet map}
\begin{align}\label{NDmap}
	 M (\lambda) \partial_\nu f_\lambda|_{\dR^n} = 
	 f_\lambda |_{\dR^n}, \quad f_\lambda \in \ker (T - \lambda), 
\end{align}
and satisfies $M (\lambda) \in \cB(L^2(\dR^n))$. 
\end{prop} 

\begin{proof}
We verify the conditions (a)--(c) of Proposition~\ref{suff_cond_qbt}. The mapping
\begin{align*}
	 H^2(\dR^{n+1}_+) \ni f \mapsto 
	\binom{\partial_\nu f|_{\dR^n}}{f |_{\dR^n}} \in 
	H^{1/2}(\dR^n) \times H^{3/2}(\dR^n)
\end{align*}
is surjective, see above. Since it is a restriction of the mapping $\Gamma = \binom{\Gamma_0}{\Gamma_1}$, the density of $H^{1/2}(\dR^n) \times H^{3/2}(\dR^n)$ in $L^2(\dR^n) \times L^2(\dR^n)$ yields (a). Condition~(b) is just the usual second Green identity,
\begin{align*}
	 \bigl(-\Delta f, g\bigr) - \bigl(f, - \Delta g\bigr) = 
	\left(f |_{\dR^n}, \partial_\nu g|_{\dR^n} \right) - 
	\left(\partial_\nu f|_{\dR^n}, g |_{\dR^n} \right),
\end{align*}
for $f, g \in H^{3/2}_\Delta (\dR_+^{n+1})$, which can be found in, e.g.,~\cite[Theorem 5.5]{Fr67}; here the inner products in $L^2 (\dR^{n + 1}_+)$ and in $L^2 (\dR^n)$ both are denoted by $(\cdot, \cdot)$. In order to verify (c) we observe that the operator $T \upharpoonright \ker \Gamma_0$ is $-\Delta$ on the domain 
\begin{align*}
 \left\{f \in H_{\Delta}^{3/2} (\dR_+^{n+1}) : \partial_\nu f|_{\dR^n} = 0 \right\} 
\end{align*}
in $L^2(\dR^{n+1}_+)$, which contains the domain of the selfadjoint Neumann operator $\AN$ in~\eqref{eq:DirNeu}. Thus by Proposition~\ref{suff_cond_qbt} $\{ L^2(\dR^n), \Gamma_0, \Gamma_1 \}$ is a quasi-boundary triple for $A^*$ and the statements on $A$ are true. In particular, $T \upharpoonright \ker \Gamma_0$ coincides with the Neumann operator $A_{\rm N}$. The representations~\eqref{Poisson} and~\eqref{NDmap} follow immediately from~\eqref{NDpropertyAbstr} and the definition of the boundary mappings $\Gamma_0$ and $\Gamma_1$. It remains to show that $M (\lambda)$ is bounded and everywhere defined. Since by Proposition~\ref{gammaprop}~(iii) $M (\lambda) \subset M (\ov \lambda)^*$ holds for each $\lambda \in \rho (\AN)$ and the latter operator is closed, $M (\lambda)$ is closable. It follows from $\dom M (\lambda) = \ran \Gamma_0 = L^2 (\dR^n)$, see~\eqref{NTrace}, that $M (\lambda)$ is even closed and, hence, $M (\lambda) \in \cB \bigl(L^2(\dR^n)\bigr)$ by the closed graph theorem.
\end{proof}

For the sake of completeness we remark that the adjoint of $A$ is given by 
\begin{align*}
 A^* f = - \Delta f, \quad \dom A^* = \left\{ f \in L^2 (\dR^{n + 1}_+) :  \Delta f \in L^2 (\dR^{n + 1}_+) \right\},
\end{align*}
but this will not play a role in our further considerations.

We are now able to provide some information on the operator $A_\alpha$ in~\eqref{eq:dom}.

\begin{thm}\label{thm1}
Let $\alpha \in W^{1, \infty} (\dR^n)$. Then each $\lambda < - \| \alpha \|_\infty^2$ belongs to $\rho (\Lar)$. Moreover, $\Lar$ is selfadjoint if and only if $\alpha$ is real-valued. In particular, in this case $\Lar$ is semibounded from below by $- \| \alpha\|_\infty^2$.
\end{thm}

\begin{rem}
We emphasize that in certain cases the estimate for the spectrum of $A_\alpha$ given in Theorem~\ref{thm1} is very rough. For example, if $\alpha$ is a real, nonpositive function, the first Green identity implies that $A_\alpha$ is even nonnegative.
\end{rem}

\begin{proof}[Proof of Theorem~\ref{thm1}]
Let $\{L^2 (\dR^n), \Gamma_0, \Gamma_1\}$ be the quasi-boundary triple for $A^*$ in Proposition~\ref{qbt}, $\gamma$ the corresponding $\gamma$-field, and $M$ the corresponding Weyl function. We verify first that with respect to the quasi-boundary triple $\{L^2 (\dR^n), \Gamma_0, \Gamma_1\}$ in Proposition~\ref{qbt} the operator $\Lar$ admits a representation $\Lar = A_{\Theta}$ in the sense of~\eqref{athetaop} with
\begin{align}\label{alphaRelation}
 \Theta = \left\{ \binom{\alpha f}{f} : f \in L^2 (\dR^n) \right\}.
\end{align}
In fact, it is obvious from the definitions that $A_\alpha \subset A_\Theta$ holds, and it remains to show $\dom A_\Theta \subset H^2 (\dR^{n + 1}_+)$. Let $f \in \dom A_\Theta \subset \dom T$ and $\eta \in \rho (A_{\rm N})$. By~\eqref{decompo} there exist $f_{\rm N}\in \dom A_{\rm N}$ and $f_\eta \in \ker (T-\eta)$ with $f = f_{\rm N} + f_\eta$. Clearly, $f_{\rm N}$ belongs to $H^2 (\dR^{n + 1}_+)$. Moreover, $f$ satisfies the boundary condition
\begin{align*}
 \alpha \Gamma_1 f = \Gamma_0 f = \Gamma_0 f_\eta;
\end{align*}
in particular, $\ran \Gamma_1 = H^1 (\dR^n)$ and the regularity assumption on $\alpha$ imply $\Gamma_0 f_\eta \in H^1 (\dR^n) \subset H^{1/2} (\dR^n)$. Since the mapping $f \mapsto \partial_\nu f |_{\dR^n}$ provides a bijection between $H^2 (\dR^{n + 1}_+) \cap \ker (T - \eta)$ and $H^{1/2} (\dR^n)$, see, e.g.,~\cite[Section~3]{G08}, it follows $f_\eta \in H^2 (\dR^{n + 1}_+)$. This shows $f \in H^2 (\dR^{n + 1}_+)$ and, hence, $A_\Theta = A_\alpha$.

Let $\lambda < - \| \alpha\|_\infty^2$ be fixed. Then $\lambda \in \rho (A_{\rm N})$ holds and by Proposition~\ref{prop.krein} in order to verify $\lambda \in \rho (A_\alpha)$ it is sufficient to show $0 \in \rho (\Theta - M(\lambda))$. Note first that $\Theta$ is injective; hence we can write
\begin{align}\label{factor}
 \left(\Theta - M (\lambda) \right)^{-1} = \Theta^{-1} \big( I - M(\lambda) \Theta^{-1} \big)^{-1} ,
\end{align}
where the equality has first to be understood in the sense of linear relations. Since $\Theta^{-1} = \alpha \in \cB(L^2(\dR^n))$, we only need to show that $I - M(\lambda) \Theta^{-1}$ has a bounded, everywhere defined inverse. In fact, the Neumann-to-Dirichlet map is given by
\begin{align*}
 M(\lambda) = (- \Delta_{\dR^n} - \lambda)^{-1/2},
\end{align*}
where $\Delta_{\dR^n}$ denotes the Laplacian in $L^2 (\dR^n)$, defined on $H^2 (\dR^n)$; cf., e.g.,~\cite[Chapter~9]{G09}. In particular, $\| M(\lambda) \| = 1 / \sqrt{- \lambda}$ holds. This implies $\| M(\lambda) \Theta^{-1} \| < 1$. Now~\eqref{factor} yields that $(\Theta - M(\lambda))^{-1}$ is a bounded operator, which is everywhere defined. This implies $\lambda \in \rho (\Lar)$. 

It follows immediately from the Green identity that $\Lar$ is symmetric if and only if $\alpha$ is real-valued. In this case $\Lar$ is even selfadjoint as $\rho (\Lar) \cap \dR$ is nonempty. Since we have shown that each $\lambda < - \| \alpha\|_\infty^2$ belongs to $\rho (\Lar)$, the statement on the semiboundedness of $A_\alpha$ follows immediately. This completes the proof.
\end{proof}

\section{Compactness and Schatten-von Neumann estimates for resolvent differences of Robin Laplacians}
\label{sec:main}

The present section is devoted to our main results on compactness and Schatten-von Neumann properties of the resolvent difference 
\begin{align}\label{resDiffagain}
 (A_{\alpha_2} - \lambda)^{-1} - (A_{\alpha_1} - \lambda)^{-1}
\end{align}
of two Robin Laplacians as in~\eqref{eq:dom} with boundary coefficients $\alpha_1$ and $\alpha_2$ in dependence of the asymptotic behavior of $\alpha_2 - \alpha_1$. Let us shortly recall the definition of the Schatten-von Neumann classes and some of their basic properties. For more details see~\cite[Chapter~II and~III]{GK69} and~\cite{S05}. Let $\sS_\infty (\cG, \cH)$ denote the linear space of all compact linear operators mapping the Hilbert space  $\cG$ into the Hilbert space $\cH$. Usually the spaces $\cG$ and $\cH$ are clear from the context and we simply write $\sS_\infty$. For $K \in \sS_\infty$ we denote by $s_k (K)$, $k = 1, 2, \dots$, the {\em singular values} (or {\em $s$-numbers}) of $K$, i.e., the eigenvalues of the compact, selfadjoint, nonnegative operator $(K^* K)^{1/2}$, enumerated in decreasing order and counted according to their multiplicities. Note that for a selfadjoint, nonnegative operator $K \in \sS_\infty$ the singular values are precisely the eigenvalues of $K$. 

\begin{defn}
An operator $K \in \sS_\infty$ is said to belong to the {\em Schatten-von Neumann class} $\sS_p$ of order $p > 0$, if its singular values satisfy
\begin{align*}
 \sum_{k=1}^\infty \bigl(s_k(K)\bigr)^p < \infty.
\end{align*}
An operator $K$ is said to belong to the {\em weak Schatten-von Neumann class} $\sS_{p, \infty}$ of order $p > 0$, if 
\begin{align*}
 s_k (K) = O (k^{- 1/p}), \quad k \to \infty,
\end{align*}
holds.
\end{defn}

Some well-known properties of the Schatten-von Neumann classes are collected in the following lemma; for proofs see the above-mentioned references and~\cite[Lemma~2.3]{BLL10}.

\begin{lem}
\label{splemma}
For $p, q, r > 0$ the following assertions hold.
\begin{itemize}\setlength{\itemsep}{1.2ex}
 \item[(i)] Let $\frac{1}{p} + \frac{1}{q} = \frac{1}{r}$. If $K \in \sS_p$ and $L \in \sS_q$, then $K L \in \sS_r$; if $K \in \sS_{p, \infty}$ and $L \in \sS_{q, \infty}$, then $K L \in \sS_{r, \infty}$;
 \item[(ii)] $K \in \sS_p~\Longleftrightarrow~K^* \in \sS_p$ and $K \in \sS_{p,\infty}~\Longleftrightarrow~K^* \in \sS_{p,\infty}$;
 \item[(iii)] $\sS_p \subset \sS_{p,\infty}$ and $\sS_{p, \infty} \subset \sS_q$ for all $q > p$, but  $\sS_{p, \infty} \not \subset \sS_p$.
\end{itemize}
\end{lem}

Let us now come to the investigation of compactness and Schatten-von Neumann properties of~\eqref{resDiffagain}. The condition 
\begin{align}\label{eq:alpha}
 \mu \left( \left\{ x \in \dR^n : | \alpha (x) | \geq \eps \right\} \right) < \infty \quad \text{for all}~\eps > 0
\end{align}
for $\alpha = \alpha_2 - \alpha_1$ turns out to be sufficient for the compactness of the resolvent difference~\eqref{resDiffagain}, see Theorem~\ref{thm2} below; here $\mu$ denotes the Lebesgue measure on $\dR^n$. We remark that the condition~\eqref{eq:alpha} includes, e.g., the case that $\alpha$ belongs to $L^q (\dR^n)$ for some $q \geq 1$, and the case that $\sup_{|x| \geq r} | \alpha (x)| \to 0$ as $r \to \infty$.

The following lemma contains the main ingredients of the proof of Theorem~\ref{thm2} and Theorem~\ref{thm2'} below. 

\begin{lem}
\label{lem:ideal}
Let $\cK$ be a Hilbert space and let $K\in\cB(\cK,L^2(\dR^n))$ be an operator with $\ran K \subset H^{3/2} (\dR^n)$. Assume $\alpha\in L^\infty(\dR^n)$.
\begin{itemize}
\item[(i)] If $\alpha$ satisfies the condition~\eqref{eq:alpha}, then $\alpha K \in \sS_{\infty}$.
\item[(ii)] If  $\alpha$ has a compact support or if $n > 3$ and $\alpha \in L^{2n/3} (\dR^n)$, then
\begin{equation*}
 \alpha K \in \sS_{\frac{2 n}{3},\infty}.
\end{equation*}
\item[(iii)] If  $\alpha \in L^2 (\dR^n)$ and $n \geq 3$, then 
\begin{equation*}
 \alpha K \in \sS_r \quad \text{for all}~r > \frac{2 n}{3}. 
\end{equation*}
\item[(iv)] If $\alpha \in L^p (\dR^n)$ for $p \geq 2$ and $p > \frac{2 n}{3}$,
then
\begin{equation*}
 \alpha K \in \sS_{p}. 
\end{equation*}
\end{itemize}
\end{lem}

\begin{proof}
Assume first that $\alpha$ satisfies~\eqref{eq:alpha}. Then there exists a sequence $\Omega_1 \subset \Omega_2 \subset \dots$ of smooth domains of finite measure whose union is all of $\dR^n$ such that for each
$m\in\dN$ we have $|\alpha(x)| < \tfrac{1}{m}$ for all $x\in \dR^n\setminus\Omega_m$. 
For each $m \in \dN$ let $\chi_m$ be the characteristic function of the set $\Omega_m$. Denote by $P_m$ the canonical projection from $L^2 (\dR^n)$ to $L^2 (\Omega_{m})$ and by $J_m$ the canonical embedding of $L^2 (\Omega_{m})$ into $L^2 (\dR^n)$. Then $\ran (P_m \chi_m K) \subset H^{3/2} (\Omega_m) \subset H^1 (\Omega_{m})$ and, by embedding statements, $P_m \chi_m K : \cK \to L^2 (\Omega_{m})$ is compact; see~\cite[Theorem~3.4 and Theorem~4.11]{EE75} and~\cite[Chapter~V]{EE87}. Since $\alpha J_m$ is bounded, it turns out that $\alpha \chi_m K = \alpha J_m P_m \chi_m K$ is compact. From the assumption~\eqref{eq:alpha} on $\alpha$ it follows easily that the sequence of operators $\alpha \chi_m K$ converges to $\alpha K$ in the operator-norm topology. Thus also $\alpha K$ is compact, which is the assertion of item~(i).

Let us assume that $\alpha$ has a compact support and that $\Omega\subset \dR^n$ is a bounded, smooth domain with $\Omega \supset \supp \alpha$. Let $P$ be the canonical projection in $L^2(\dR^n)$ onto $L^2(\Omega)$ and let $J$ be the canonical embedding of $L^2(\Omega)$ into  $L^2(\dR^n)$, and let $\wt \alpha := \alpha |_{\Omega}$. Since $\ran(PK) \subset H^{3/2}(\Omega)$ and $\Omega$ is a bounded, smooth domain, the embedding operator from $H^{3/2} (\Omega)$ into $L^2(\Omega)$ is contained in the class $\sS_{\frac{2 n}{3},\infty}$, see~\cite[Theorem 7.8]{HT03}. It follows $PK\in \sS_{\frac{2 n}{3},\infty}$ as a mapping from $\cK$ into $L^2(\Omega)$. Since $J \wt \alpha$ is bounded, we obtain $\alpha K = J \wt \alpha P K \in \sS_{\frac{2 n}{3},\infty}$.

The proofs of the remaining statements make use of spectral estimates for the operator $\alpha D$ in $L^2 (\dR^n)$ with
\begin{align}\label{g}
 D = (I - \Delta_{\dR^n})^{- 3/4} = g (- i \nabla), \quad g (x) = (1 + |x|^2)^{- 3/4},~ x \in \dR^n,
\end{align}
where the formal notation $g (- i \nabla)$ can be made precise with the help of the Fourier transformation. We remark that $D$, regarded as an operator from $L^2 (\dR^n)$ into $H^{3/2} (\dR^n)$, is an isometric isomorphism. Recall that a function $f$ is said to belong to the {\em weak Lebesgue space} $L^{p, \infty} (\dR^n)$ for some $p > 0$, if the condition
\begin{align*}
 \sup_{t > 0} \big( t^p \mu \big(\{x \in \dR^n : |f(x)| > t \} \big) \big) < \infty
\end{align*}
is satisfied, where $\mu$ denotes the Lebesgue measure on $\dR^n$. The function $g$ in~\eqref{g} belongs to $L^{2 n/3,\infty}(\dR^n)$. In fact, one easily verifies that the set $\{x \in \dR^n : | g (x)| > t\}$ is contained in the ball of radius $t^{- 2/3}$ centered at the origin, and the formula for the volume of a ball leads to the claim. Let now $ n> 3$ and $\alpha \in L^{2n/3} (\dR^n)$. Then a result by M.\,Cwikel in~\cite{C77} yields
\begin{align*}
 \alpha D \in \sS_{\frac{2 n}{3}, \infty};
\end{align*}
see also~\cite[Theorem~4.2]{S05}. We conclude
\begin{align*}
 \alpha K = \alpha D D^{-1} K \in \sS_{\frac{2 n}{3}, \infty}.
\end{align*}
Thus we have proved~(ii).

In order to show~(iii) let us assume $\alpha \in L^2 (\dR^n)$ and $n \geq 3$.
Since $\alpha$ is bounded,  $\alpha\in L^p(\dR^n)$ for each $p > 2$.
It is easy to check that $g$ in~\eqref{g} belongs to $L^p (\dR^n)$ for each $p > 2 n/3$.  The standard result~\cite[Theorem~4.1]{S05} and  $\alpha, g  \in L^r (\dR^n)$ for all $r > 2n/3\ge 2$ imply
\begin{align*}
 \alpha D \in \sS_r, \quad \text{for all}~ r > \frac{2 n}{3}.
\end{align*}
It follows
\begin{align*}
 \alpha K = \alpha D D^{-1} K \in \sS_r, \quad \text{for all}~ r > \frac{2 n}{3},
\end{align*}
which is the assertion of~(iii).

Let now $\alpha \in L^p (\dR^n)$ for $p \geq 2$ and $p > 2 n/3$. As above, $g \in L^p (\dR^n)$ and~\cite[Theorem~4.1]{S05} yields $\alpha D \in \sS_p$. Hence, $\alpha K = \alpha D D^{-1} K \in \sS_p$, which completes the proof of~(iv).
\end{proof}

\begin{rem}
The condition in Lemma~\ref{lem:ideal}~(i) can still be slightly weakened using the optimal prerequisites on a domain $\Omega$ which imply compactness of the embedding of $H^1 (\Omega)$ into $L^2 (\Omega)$; see, e.g.,~\cite[Chapter~VIII]{EE87}. To avoid too inconvenient and technical assumptions, we restrict ourselves to the above condition.
\end{rem}

We continue with giving a factorization of the resolvent difference of two Robin Laplacians. It is based on the formula of Krein type in Proposition~\ref{prop.krein} and will be crucial for the proofs of our main results. We remark that an analogous formula as below is well known for ordinary boundary triples and abstract boundary conditions, see~\cite[Proof~of~Theorem~2]{DM91}.

\begin{lem}\label{fact}
Let $\alpha_1, \alpha_2 \in W^{1, \infty} (\dR^n)$ and let $A_{\alpha_1}, A_{\alpha_2}$ be the corresponding Robin Laplacians as in~\eqref{eq:dom}. Then
\begin{align*}
 (A_{\alpha_2} - \lambda & )^{-1} - (A_{\alpha_1} - \lambda)^{-1}\\
 & = \gamma (\lambda) \left( I - \alpha_1 M (\lambda) \right)^{-1} \left( \alpha_2 - \alpha_1 \right) \left( I - M (\lambda) \alpha_2 \right)^{-1} \gamma ( \lambda)^*
\end{align*}
holds for each $\lambda < - \max \{ \|\alpha_1\|_\infty^2, \|\alpha_2\|_\infty^2 \}$, where $\gamma (\lambda)$ is the Poisson operator in~\eqref{Poisson} and $M (\lambda)$ is the Neumann-to-Dirichlet map in~\eqref{NDmap}.
\end{lem}

\begin{proof}
Let $A$ be given as in~\eqref{eq:A} and let $\{L^2 (\dR^n), \Gamma_0, \Gamma_1\}$ be the quasi-boundary triple for $A^*$ in Proposition~\ref{qbt}, so that $\gamma$ is the corresponding $\gamma$-field and $M$ is the corresponding Weyl function. Let us fix $\lambda$ as in the proposition. Then $\lambda$ belongs to $\rho (A_{\alpha_1}) \cap \rho (A_{\alpha_2})$ by Theorem~\ref{thm1}. Moreover, if $\Theta_1$ and $\Theta_2$ denote the  linear relations corresponding to $\alpha_1$ and $\alpha_2$, respectively, as in~\eqref{alphaRelation}, then we have
\begin{align*}
 \big(\Theta_2 & - M (\lambda) \big)^{-1} - \big(\Theta_1 - M (\lambda) \big)^{-1}\\ 
 & = \alpha_2 \big(I - M (\lambda) \alpha_2 \big)^{-1} - \big( I - \alpha_1 M (\lambda) \big)^{-1} \alpha_1 \\
 & = \big( I - \alpha_1 M (\lambda) \big)^{-1} \Big(\big( I - \alpha_1 M (\lambda) \big) \alpha_2 - \alpha_1 \big(I - M (\lambda) \alpha_2 \big) \Big) \big(I - M (\lambda) \alpha_2 \big)^{-1},
\end{align*}
which, together with Proposition~\ref{prop.krein}, completes the proof.
\end{proof}

The following two theorems contain the main results of the present paper. Since their proofs have similar structures, we give a joint proof below. The first of the two main theorems states that under the condition~\eqref{eq:alpha} on $\alpha = \alpha_2 - \alpha_1$ the resolvent difference~\eqref{resDiffagain} is compact.

\begin{thm}\label{thm2}
Let $\alpha_1, \alpha_2 \in W^{1, \infty} (\dR^n)$, let $A_{\alpha_1}$, $A_{\alpha_2}$ be the corresponding operators as in~\eqref{eq:dom}, and let $\alpha := \alpha_2 - \alpha_1$ satisfy~\eqref{eq:alpha}. Then 
 \begin{align*}
  (A_{\alpha_2} - \lambda)^{-1} - (A_{\alpha_1} - \lambda)^{-1} \in \sS_\infty
 \end{align*}
holds for each $\lambda \in \rho (A_{\alpha_1}) \cap \rho (A_{\alpha_2})$, and, in particular, $\sigma_{\rm ess} (A_{\alpha_1}) = \sigma_{\rm ess} (A_{\alpha_2})$.
\end{thm}

As mentioned before, the condition~\eqref{eq:alpha} covers the case that $\alpha$ belongs to $L^p (\dR^n)$ for an arbitrary $p > 0$. For certain $p$, if $\alpha \in L^p (\dR^n)$, the result of Theorem~\ref{thm2} can be improved as follows.

\begin{thm}\label{thm2'}
Let $\alpha_1, \alpha_2 \in W^{1, \infty} (\dR^n)$, let $A_{\alpha_1}$, $A_{\alpha_2}$ be the corresponding operators as in~\eqref{eq:dom}, and let $\alpha := \alpha_2 - \alpha_1$. Then for $\lambda \in \rho (A_{\alpha_1}) \cap \rho (A_{\alpha_2})$ the following assertions hold.
\begin{enumerate}
 \item If $\alpha$ has a compact support or if $n> 3$ and $\alpha \in L^{n/3} (\dR^n)$, then
 \begin{equation*}
  (A_{\alpha_2} - \lambda)^{-1}-(A_{\alpha_1} -\lambda)^{-1}\in\sS_{\frac{n}{3},\infty}.
\end{equation*}
 \item If $\alpha \in L^1 (\dR^n)$ and $n \geq 3$, then
 \begin{align*}
  (A_{\alpha_2} - \lambda)^{-1}-(A_{\alpha_1} -\lambda)^{-1}\in\sS_{r} \quad \text{for all}~r > \frac{n}{3}.
 \end{align*}
 \item If $\alpha \in L^p (\dR^n)$ for $p \geq 1$ and $p > n/3$, then
 \begin{align*}
  (A_{\alpha_2} - \lambda)^{-1}-(A_{\alpha_1} -\lambda)^{-1}\in\sS_{p}.
 \end{align*}
\end{enumerate}
\end{thm}

\begin{proof}[Proof of Theorem~\ref{thm2} and Theorem~\ref{thm2'}]
Let us fix $\lambda < - \max \{ \|\alpha_1\|_\infty^2, \|\alpha_2\|_\infty^2 \}$. We first observe that 
\begin{align}\label{range}
 \ran \bigl((I - M (\lambda) \alpha_2)^{-1} \gamma ( \lambda)^*\bigr) \subset H^{3/2} (\dR^n).
\end{align}
Note first that $\dom A_{\rm N} \subset H^2 (\dR^{n + 1}_+)$ and Proposition~\ref{gammaprop} (ii) imply $\ran \gamma ( \lambda)^* \subset H^{3/2} (\dR^n)$. Thus, for $\varphi \in \ran (I - M (\lambda) \alpha_2)^{-1} \gamma (\lambda)^*$ we have $\varphi - M (\lambda) \alpha_2 \varphi \in H^{3/2} (\dR^n)$, and $\ran M (\lambda) \subset H^1 (\dR^n)$ implies $\varphi \in H^{1} (\dR^n)$. Since $\alpha_2 \varphi$ belongs to $H^1 (\dR^n)$, $M (\lambda) \alpha_2 \varphi$ automatically belongs to $H^2 (\dR^n)$; this can be seen as in the proof of Theorem~\ref{thm1}, see also~\cite[Section~3]{G08}. This proves~\eqref{range}. Analogously also
\begin{align}\label{range2}
\ran \bigl((I - M (\lambda) \overline{\alpha_1})^{-1} \gamma ( \lambda)^*\bigr) \subset H^{3 /2} (\dR^n)
\end{align}
holds. The factorization given in Lemma~\ref{fact} can be written as
\begin{align}\label{factorSupport}
 (A_{\alpha_2} - \lambda & )^{-1} - (A_{\alpha_1} - \lambda)^{-1} \nonumber\\
 & = \gamma (\lambda) \left( I - \alpha_1 M (\lambda) \right)^{-1} \sqrt{|\alpha|} \widetilde \alpha \sqrt{|\alpha|} \left( I - M (\lambda) \alpha_2 \right)^{-1} \gamma ( \lambda)^*,
\end{align}
where $\widetilde \alpha (x)$ is given by $0$ if $\alpha (x) = 0$ and by $\alpha (x) / |\alpha (x)|$ if $\alpha (x) \neq 0$. 

If $\alpha$ satisfies~\eqref{eq:alpha}, then the same holds for $\alpha$ replaced by $\sqrt{|\alpha|}$. Now~\eqref{range} and Lemma~\ref{lem:ideal} (i) imply $$\sqrt{|\alpha|} \left( I - M (\lambda) \alpha_2 \right)^{-1} \gamma (\lambda)^* \in \sS_\infty.$$ Since $\gamma (\lambda) \left( I - \alpha_1 M (\lambda) \right)^{-1} \sqrt{|\alpha|} \widetilde \alpha \in \cB (L^2 (\dR^n))$, the assertion of Theorem~\ref{thm2} follows from~\eqref{factorSupport}.

If $\alpha$ has a compact support or if $ n>3$ and $\alpha$ belongs to $L^{n/3} (\dR^n)$, then $\sqrt{|\alpha|}$ has a compact support or belongs to $L^{2n/3} (\dR^n)$, respectively; thus~\eqref{range} and~\eqref{range2} together with Lemma~\ref{lem:ideal} (ii) imply 
\begin{align*}
 \sqrt{|\alpha|} (I - M (\lambda) \alpha_2)^{-1} \gamma (\lambda)^* \in \sS_{\frac{2 n}{3},\infty} \quad \text{and} \quad \sqrt{|\alpha|} (I - M (\lambda) \overline{\alpha_1})^{-1} \gamma ( \lambda)^* \in \sS_{\frac{2 n}{3},\infty}.
\end{align*}
Taking the adjoint of the latter operator, Lemma~\eqref{splemma} (i) and (ii) and~\eqref{factorSupport} yield~Theorem~\ref{thm2'}~(i).

The proofs of Theorem~\ref{thm2'}~(ii) and~(iii) are completely analogous; one uses Lemma~\ref{lem:ideal}~(iii) and~(iv), respectively, instead of item~(ii).
\end{proof}

As an immediate consequence of Theorem~\ref{thm2'} we obtain the following result concerning scattering theory. Note that in the case $n < 3$ and $\alpha_2 - \alpha_1 \in L^1 (\dR^n)$ Theorem~\ref{thm2'}~(iii) implies that the difference~\eqref{resDiffagain} is contained in the trace class $\sS_1$. Now basic statements from scattering theory yield the following corollary, see, e.g.,~\cite[Theorem X.4.12]{K95}.

\begin{cor}
Let $n < 3$ and let $\alpha_1, \alpha_2\in W^{1,\infty}(\dR^n)$ be real-valued with $\alpha_2 - \alpha_1 \in L^1(\dR^n)$. Then wave operators for the pair of selfadjoint operators $\{A_{\alpha_1}, A_{\alpha_2}\}$ exist and are complete. Moreover, the absolutely continuous spectra of $A_{\alpha_1}$ and $A_{\alpha_2}$ coincide and their absolutely continuous parts are unitarily equivalent.
\end{cor}

We would like to put some emphasis on the important special case $\alpha_1 = 0$, in which $A_{\alpha_1}$ is the selfadjoint Neumann operator $A_{\rm N}$ in~\eqref{eq:DirNeu}. In this situation Theorem~\ref{thm2} and Theorem~\ref{thm2'} read as follows. Recall that the spectrum of $A_{\rm N}$ has the simple structure $\sigma (A_{\rm N}) = \sigma_{\rm ess} (A_{\rm N}) = [0, \infty)$.

\begin{cor}\label{thmN}
Let $\alpha \in W^{1, \infty} (\dR^n)$ satisfy~\eqref{eq:alpha} and let $A_{\alpha}$ be the operator in~\eqref{eq:dom}. Then 
 \begin{align*}
  (A_{\alpha} - \lambda)^{-1} - (A_{\rm N} - \lambda)^{-1} \in \sS_\infty
 \end{align*}
holds for each $\lambda \in \rho (A_{\alpha_1}) \cap \rho (A_{\alpha_2})$, and, in particular, $\sigma_{\rm ess} (A_{\alpha}) = [0, \infty)$.
\end{cor}

\begin{cor}\label{thmN'}
Let $\alpha \in W^{1, \infty} (\dR^n)$ and let $A_{\alpha}$ be the operator in~\eqref{eq:dom}. Then for $\lambda \in \rho (A_{\alpha}) \cap \rho (A_{\rm N})$ the following assertions hold.
\begin{enumerate}
 \item If $\alpha$ has a compact support or if $n > 3$ and $\alpha \in L^{n/3} (\dR^n)$, then
 \begin{equation*}
  (A_{\alpha} - \lambda)^{-1}-(A_{\rm N} -\lambda)^{-1}\in\sS_{\frac{n}{3},\infty}.
\end{equation*}
 \item If $\alpha \in L^1 (\dR^n)$ and $n \geq 3$, then
 \begin{align*}
  (A_{\alpha} - \lambda)^{-1}-(A_{\rm N} -\lambda)^{-1}\in\sS_{r} \quad \text{for all}~r > \frac{n}{3}.
 \end{align*}
 \item If $\alpha \in L^p (\dR^n)$ for $p \geq 1$ and $p > n/3$, then
 \begin{align*}
  (A_{\alpha} - \lambda)^{-1}-(A_{\rm N} -\lambda)^{-1}\in\sS_{p}.
 \end{align*}
\end{enumerate}
\end{cor}

As a consequence of Corollary~\ref{thmN}, applying~\cite[Proposition~5.11~(v) and~(vii)]{MN11} we obtain the following statement on the absolutely continuous part of $A_\alpha$, if $\alpha$ is real-valued.

\begin{cor}
Let $\alpha \in W^{1, \infty} (\dR^n)$ be real-valued satisfying~\eqref{eq:alpha} and let $A_{\alpha}$ be the selfadjoint operator in~\eqref{eq:dom}. Then $A_{\rm N}$ and the absolutely continuous part of $A_\alpha$ are unitarily equivalent.
\end{cor}

We conclude our paper with an observation connected with the speed of accumulation of the discrete spectrum of the operator $A_\alpha$, where $\alpha$ is a complex-valued function subject to the condition~\eqref{eq:alpha}. As Corollary~\ref{thmN} shows, the essential spectrum of $A_\alpha$ in this case is given by $[0, \infty)$ and, additionally, discrete, (in general) non-real eigenvalues may appear. The following statement combines our main result with some recent advances in the theory of non-selfadjoint perturbations of selfadjoint operators; it is based on~\cite[Theorem~2.1]{H11}.

\begin{cor}
Let $\alpha \in W^{1, \infty} (\dR^n)$ and let $A_{\alpha}$ be the operator in~\eqref{eq:dom}. Then for all $a > \|\alpha\|_\infty^2$ the following assertions hold.
\begin{enumerate}
\item If $\alpha \in L^p (\dR^n)$ for $p \geq 1$ and $p > n/3$, then
\[
 \sum_{\lambda\in\sigma_{\rm d}(A_\alpha)}\dist\bigl((\lambda+a)^{-1}, [0,a^{-1}]\bigr)^p < \infty.
\]
\item If $n\ge 3$ and $\alpha \in L^{n/3} (\dR^n)$, then
\[
 \sum_{\lambda\in\sigma_{\rm d}(A_\alpha)} \dist\bigl((\lambda+a)^{-1}, [0,a^{-1}]\bigr)^{\tfrac{n}{3}+\varepsilon} < \infty \quad \text{for all}\quad \varepsilon > 0.
\]
\end{enumerate}
Above the eigenvalues in the discrete spectrum are counted according to their algebraic multiplicities, and $\dist(\cdot,\cdot)$ denotes the usual distance in the complex plane.
\end{cor}

For the proof recall that the numerical range of a bounded operator $A$ in a Hilbert space $\cH$ is defined as
\[
{\rm Num}(A) := \left\{ (Af,f)_{\cH}\colon f\in\cH, \|f\|_\cH = 1\right\}.
\]

\begin{proof}
Let us assume first $\alpha \in L^p (\dR^n)$ for some $p \geq 1$ with $p > n /3$. Clearly $-a\in\rho(A_{\rm N})$ and by Theorem~\ref{thm1} also $-a\in\rho(A_\alpha)$. In view of the assumptions on $\alpha$ it follows from Corollary~\ref{thmN'}~(iii) that
\begin{equation}
\label{Sp}
(a + A_\alpha)^{-1} - (a+ A_{\rm N})^{-1} \in \sS_p. 
\end{equation}
The operator $(a + A_{\rm N})^{-1}$ is bounded and selfadjoint and it has a purely essential spectrum given by $\big[0,a^{-1}\big]$. Since
\begin{equation*}
\ov{{\rm Num}\big((a + A_{\rm N})^{-1}\big)} = \sigma \left( (a + A_{\rm N})^{-1} \right) = \big[0,a^{-1}\big]
\end{equation*}
and, trivially,
\begin{equation*}
\lambda\in\sigma_{\rm d}(A_\alpha) \quad \Longleftrightarrow \quad (\lambda + a)^{-1} \in \sigma_{\rm d}\big((a+A_\alpha)^{-1}\big),
\end{equation*}
the claim of~(i) follows from \cite[Theorem 2.1]{H11}. 

The proof of (ii) uses Corollary~\ref{thmN'}~(i) and~(ii) instead of item~(iii) and is completely analogous. 
\end{proof}

\end{document}